\documentclass[a4paper,10pt]{article}

\usepackage{amssymb}

\newtheorem{thm}{Theorem}[section]
\newtheorem{cor}{Corollary}[section]
\newtheorem{lem}{Lemma}[section]
\newtheorem{prop}{Proposition}[section]

\newtheorem{exmp}{Example}[section]
\newtheorem{rem}{Remark}[section]

\newcommand{\qed}{\square}
\renewcommand{\(}{\left(}
\renewcommand{\)}{\right)}

\newcommand{\im}{{\rm im}}

\newcommand{\K}{\wh{K}}
\newcommand{\got}{\mathfrak}

\newcommand{\X}{X_1,\ldots ,X_n}

\newcommand{\Ym}{Y_1,\ldots ,Y_m}
\newcommand{\Y}{Y_1,\ldots ,Y_n}

\renewcommand{\v}{\wh{v}}

\newcommand{\wh}{\widehat}
\newcommand{\om}{\omega}
\newcommand{\Om}{\Omega}
\newcommand{\al}{\alpha}

\newcommand{\si}{\sigma}

\newcommand{\De}{\Delta}

\newcommand{\be}{\beta}
\newcommand{\vp}{\varphi}

\newcommand{\R}{R_{\v}}
\renewcommand{\r}{R_{v}}
\newcommand{\M}{{\got m}_{\v}}
\newcommand{\m}{{\got m}_v}

\newcommand{\Kn}{{\widehat K_n}}
\newcommand{\D}{\Delta_{\v}}
\renewcommand{\d}{\Delta_v}

\newcommand{\Phii}{\Phi^{-1}}
\newcommand{\pst}{\Phi_{\sigma ,\theta}}
\newcommand{\psti}{\(\pst\)^{-1}}
\newcommand{\lcor}{{\rm [\kern - 1.8pt [}}
\newcommand{\rcor}{{\rm ]\kern - 1.8pt ]}}
\renewcommand{\[}{\left[}
\renewcommand{\]}{\right]}
\newcommand{\llcor}{\[\kern -2.5pt\[}
\newcommand{\rrcor}{\]\kern -2.5pt\]}
\newcommand{\Dt}{\D\lcor t\rcor}

\newcommand{\fb}{\overline{f}}

\newcommand{\ZZ}{{\mathbb{Z}}}
\newcommand{\QQ}{{\mathbb{Q}}}
\newcommand{\RR}{{\mathbb{R}}}
\newcommand{\CC}{{\mathbb{C}}}
\newcommand{\FF}{{\mathbb{F}}}

\newcommand{\LL}{{\mathbb{L}}}
\newcommand{\lto}{\longrightarrow}
\newcommand{\lmapsto}{\longmapsto}
\title{Rank one discrete valuations of power series fields}
\author{F. J. Herrera Govantes\thanks{Partially supported by MTM2004-07203-C02-01 and FEDER}\\
Departamento de \'Algebra\\
Universidad de Sevilla\\
email: jherrera@us.es \and
M. A. Olalla Acosta$^*$\\
Departamento de \'Algebra\\
Universidad de Sevilla\\
email: miguelolalla@us.es \and
J.L. Vicente C\'ordoba$^*$\\
Departamento de \'Algebra\\
Universidad de Sevilla\\
email: jlvc@us.es}

\begin{document}

\maketitle

\begin{abstract}
In this paper we study rank one discrete valuations of the field\newline
$k((X_1,\ldots ,X_n))$ whose center in $k\lcor\X\rcor$ is the maximal
ideal. In sections 2 to 6 we give a construction of a
system of parametric equations describing such valuations. This amounts to
finding a parameter and a field of coefficients. We devote section 2
to finding an element of value 1, that is, a parameter. The field of
coefficients is the residue field of the valuation, and it is
given in section 5.

The constructions given in these sections are not effective in the
general case, because we need either to use Zorn's lemma or to know
explicitly a section $\sigma$ of the natural homomorphism $R_v\to\d$
between the ring and the residue field of the valuation $v$.

However, as a consequence of this construction, in section 7, we
prove that $k((\X ))$ can be embedded into a field $L((\Y ))$, where
$L$ is an algebraic extension of $k$ and the {\em ``extended valuation'' is as close as possible to the usual
  order function}.
\end{abstract}


\section{Terminology and preliminaries}

Let $k$ be a field of characteristic 0, $R_n=k\lcor\X\rcor$ the formal power series ring in $n$ variables, $M_n=(\X )$ its maximal ideal and $K_n=k ((\X ))$ its quotient field. Let $v$ be a rank-one discrete valuation of $K_n\vert k$, $\r$ the valuation ring, $\m$ the maximal ideal and $\d$ the residue field of $v$. The center of $v$ in $R_n$ is $\m\cap R_n$. Throughout this paper ``discrete valuation of $K_n\vert k$" means ``rank-one discrete valuation of $K_n\vert k$ whose center in $R_n$ is the maximal ideal $M_n$". The dimension of $v$, $\dim (v)$, is the transcendence degree of $\d$ over $k$.

For a valuation $v$ of $K\vert k$ with value group $\Gamma_v$ the well known {\em Abhyankar's inequality} \cite{Abh} says that
$$\dim (K\vert k)\ge \mbox{rat.rank}(\Gamma )_v+\dim (v),$$
where $\mbox{rat.rank}(\Gamma_v)$ is the rational rank of $v$, i.e. the dimension of $\Gamma_v\otimes\QQ$ as a $\QQ$--vector space. In our case $\dim (K_n\vert k)=n$ and $\mbox{rat.rank}(\Gamma_v)=1$, so $\dim (v)\leq n-1$. We actually know \cite{MAO} that $\dim (v)$ can be any number between 1 and $n-1$.

In order to simplify the writing we shall assume, without loss of generality, that the group of $v$ is $\ZZ$.

For $\alpha\in\ZZ$, consider the following $\r$--submodules of $K_n$:
$${\bf P}_\alpha=\{f\in K_n\mid v(f)\ge\alpha\}\cup\{0\},$$
$${\bf P}_{\alpha+}=\{f\in K_n\mid v(f)>\alpha\}\cup\{0\}.$$
We define
$$G_v=\bigoplus_{\alpha\in\ZZ}\frac{{\bf P}_\alpha}{{\bf P}_{\alpha+}}.$$
The $\d$--algebra $G_v$ is an integral domain. For any element $f\in K_n$ with $v(f)=\alpha$, the natural image of $f$ in $\frac{{\bf P}_\alpha}{{\bf P}_{\alpha+}}\subset G_v$ is a homogeneous element of $G_v$ of degree $\alpha$ called the $v$--initial form of $f$, which we will denote by $in_vf$. The algebra $G_v$ is saturated.

Let $\Kn$ be the completion of $K_n$ with respect to $v$, $\v$ the extension of $v$ to $\Kn$, $\R$, $\M$ and $\D$ the ring, maximal ideal and the residue field of $\v$, respectively (see \cite{Ser1} for more details). We know that $\d$ and $\D$ are isomorphic \cite{Krull}. Let $\si :\D\to\R$ be a $k-$section of the natural homomorphism $\R\to\D$, $\theta\in\R$ an element of value 1 and $t$ an indeterminate. Since Kaplansky \cite{Kap} we know that $\Dt$ and $\R$ are isomorphic, let us consider the $k-$isomorphism
$$
\Phi =
\pst :\Dt\to\R
$$
given by
$$
\Phi\( \sum \al_i t^i \) =\sum\si (\al_i)\theta^i,
$$
and denote also by $\Phi$ its extension to the quotient fields. Then we have a $k-$isomorphism $\Phii$ which, when composed with the usual order function on $\D ((t))$, gives the valuation $\v$. This is the situation we will consider throughout this paper, and we will freely use it without new explicit references.

We shall use two basic transformations in order to find an element of value 1 and construct the residue field:

\begin{enumerate}

\item {\em Monoidal transformation:}
$$ \begin{array}{rcl} k\lcor\X\rcor & \lto & k\lcor\Y\rcor \\
X_1 & \lmapsto & Y_1 \\
X_2 & \lmapsto & Y_1Y_2 \\
X_i & \lmapsto & Y_i,\ i=3,\ldots ,n.
\end{array}
$$ with $v(X_2)>v(X_1)$.

\item {\em Change of coordinates:}
$$\begin{array}{rcl} k\lcor\X\rcor & \lto & L\lcor\Y\rcor \\
X_1 & \lmapsto & Y_1 \\ 
X_i & \lmapsto & Y_i+c_iY_1,\ i=2,\ldots ,n,
\end{array}$$
where $c_i\in\R\setminus\M$ and $L$ is an extension field of $k$.

\end{enumerate}

For both transformations we have the following facts:

\begin{enumerate}

\item[(a)] The transformations are one to one. In the case of the monoidal transformations this property is well known. In the other case it is a consequence of \cite{ZSII} (corollary 2, page 137).


\item[(b)] The new variables $Y_i$ lie in $\R$, so we can put $\psti (Y_i)=\sum a_{i,j}t^j$.

\item[(c)] Let $\varphi :R_n\to\d ((t))$ be the restriction of $\psti$ to $K_n$. Let us denote by $\varphi ':L\lcor\Y\rcor\to\d ((t))$ the natural extension of $\varphi$ to the ring $L\lcor\Y\rcor$. We shall denote by $\phi'$ the natural extension of $\varphi '$ to its quotient field $L_n=L((\Y ))$. Then $v=\nu_t\circ\phi '_{\vert K_n}$, with $\nu_t$ the usual order function over $\d ((t))$. Therefore, if $\varphi '$ is injective we can extend the valuation $v$ to the field $L_n$ and the extension is $v'=\nu_t\circ\phi '$.

\end{enumerate}

\begin{rem}
In fact, a monoidal transformation is a blowing-up followed by the completion of a local ring with respect to its maximal ideal. This is a problem because, in the general case, the valuation of a local ring can not be extended to the completion \cite{HOST}. That is why the homomorphism $\varphi '$ could be not injective. In this case the rank one valuation can be extended to the quotient ring $L\lcor\Y\rcor /\ker (\varphi ')$. We shall see some examples later.
\end{rem}

From now on transformation will mean monoidal transformation, change of coordinates, variable permutations or finite compositions of these.


\section{Finding an element of value 1}

Remember that we are assuming that the group of $v$ is $\ZZ$, so there exists an element $u\in K$ such that $v(u)=1$.

\begin{lem}{\label{lema1}}
Let $\al_i=v(X_i)$ for all $i=1,\ldots, n$. By means of a finite number of monoidal transformations we can find $n$ elements $Y_1,\ldots ,Y_n\in\Kn$ such that $v (Y_i)=\al =\gcd\{\al_1,\ldots ,\al_n\}$.
\end{lem}

\noindent{\it Proof.}
We can suppose that $v(X_1)=\al_1=\min\{\al_i\vert 1\leq i\leq n\}$ and consider the following two steps:

{\bf Step 1.-} If there exists $n_i\in\ZZ$ such that $\al_i=n_i\al_1$ for all $i= 2,\ldots , n$, then for each $i$ we apply $n_i-1$ monoidal transformations
$$ \begin{array}{rcl} k\lcor\X\rcor & \lto &
k\lcor\Y\rcor \\
X_i & \lmapsto & Y_1Y_i\\
X_j & \lmapsto & Y_j,\ j\ne i.
\end{array}
$$
Trivially $v (Y_i)=\al_1$ for all $i=1,\ldots ,n$.

{\bf Step 2.-} Assume there exists $i$, with $2\leq i\leq n$, such that $v(X_1)=\al_1$ does not divide to $v(X_i)=\al_i$. We can suppose that $i=2$ with no loss of generality and then $\al_2 = q\al_1+r$. So we apply $q$ times the monoidal transformation
$$
\begin{array}{rcl} k\lcor\X\rcor & \lto & k\lcor\Y\rcor \\
X_2 & \lmapsto & Y_1Y_2 \\
X_i & \lmapsto & Y_i,\
i\ne 2
\end{array}
$$
to obtain a new ring $k \lcor\Y\rcor$ where $v (Y_2) =r>0$ and $Y_2$ is the element of minimum value.

As the values of the variables are greater than zero, in a finite number of steps 2 we come to the situation of step 1. In fact, this algorithm is equivalent to the ``euclidean algorithm" to compute the greatest common divisor of $\al_1,\ldots ,\al_n$. $\qed$

\begin{thm}
We can find an element of value 1 applying a finite number of monoidal transformations and changes of coordinates.
\end{thm}

\noindent{\it Proof.}
We shall denote $Y_{1,r},\ldots ,Y_{n,r}$ the elements found after $r$ transformations.

We can suppose that we have applied the previous lemma to obtain some elements $Y_{1,1},\ldots ,Y_{n,1}$ such that $v(Y_{i,1})=\al$ for all $i=1,\ldots ,n$. Let us prove that there exists $c_i\in\R\setminus\M$ for each $i=2,\ldots ,n$ such that $\v (Y_{i,1}-c_iY_{1,1}) >\al$. We can take
$$\psti (Y_{i,1})=\sum_{j\geq\al} a_{i,j}t^j=\om_i(t),\
a_{i,j} \in\D,\ a_{i,\al}\ne 0,$$
and so it suffices taking $b_i=a_{i,\al}/a_{1,\al }$ and $c_i=\si (b_i)$. In the graduated $\D$--algebra $G_{\v}$, $in_{\v}(Y_{i,1})=c_iY_{1,1}$. So the following procedure search the modules $\frac{{\bf P}_\alpha}{{\bf P}_{\alpha+}}\subset G_{\v}$ for the initial forms of the elements $Y_{i,r}$.

The following two steps define a procedure to obtain an element of value 1:

{\bf Step 1.-} We apply the coordinate change
$$
\begin{array}{rcl}
k\lcor Y_{1,1},\ldots ,Y_{n,1}\rcor & \lto & L\lcor Y_{1,2},\ldots
,Y_{n,2}\rcor \\
Y_{1,1} & \lmapsto & Y_{1,2} \\
Y_{i,1} & \lmapsto & Y_{i,2}+ c_iY_{1,2},\ i=2,\ldots ,n.
\end{array}
$$
With this transformation the values of the new variables are not equal to $\v (Y_{1,2})$.

{\bf Step 2.-} We apply lemma \ref{lema1} to obtain the same values for the elements and go to step 1. Obviously, the minimum of the values of the elements does not increase, because the greater common divisor of the values does not exceed the minimum of the values. Moreover the first variable does not change.

If we obtain an element of value 1 then we are finished.

We have to show that the procedure produces an element of value 1 in a finite number of transformations. The only way for the process to be infinite is that, in step 2, the minimum of the values of the elements does not decrease. This means that, in step 1, the value of the first variable divides the values of the new variables.

The composition of steps 1 and 2 is the transformation
$$
\begin{array}{rcl}
k\lcor Y_{1,r},\ldots ,Y_{n,r}\rcor & \lto & L\lcor Y_{1,r+1}\ldots
,Y_{n,r+1}\rcor \\
Y_{1,r} & \lmapsto & Y_{1,r+1} \\
Y_{i,r} & \lmapsto & Y_{i,r+1}+c_iY_{1,r+1}^{m_i},\ i=2,\ldots ,n.
\end{array}
$$
Where $c_iY_{1,r}^{m_i}=in_{\v}(Y_{i,r})$ in $G_{\v}$. If $\v (Y_{i,r})=\al_{i,r}$ then $c_iY_{1,r}^{m_i}$ is a homogeneous form of ${\bf P}_{\al_{i,r}}\setminus{\bf P}_{\al_{i,r}+}$ and $\v (Y_{i,r+1})=\v (Y_{i,r}-c_iY_{1,r}^{m_i})=\al_{i,r+1}>\al_{i,r}$

If we use steps 1 and 2 infinitely many times, we have an infinite sequence of transformations
$$
\begin{array}{rcl}
k\lcor\Y\rcor & \lto & L\lcor Y_{1,j},\ldots ,Y_{n,j}\rcor \\
Y_1 & \lmapsto & Y_{1,j}
\\ Y_i & \lmapsto & Y_{i,j}+\sum_{k=1}^j c_{i,k}Y_{1,j}^{m_{i,k}},\
i=2,\ldots ,n.
\end{array}
$$
Where $c_{i,j}Y_1^{m_{i,j}}=in_{\v}(Y_i-\sum_{k=1}^{j-1} c_{i,k}Y_1^{m_{i,k}})$ in $G_{\v}$. So, if $\v (Y_{i,k})=\al_{i,k}$, $\sum_{k=1}^{j} c_{i,k}Y_1^{m_{i,k}}\in {\bf P}_{\al_{i,j-1}}\setminus{\bf P}_{\al_{i,j-1}+}$ and
$$
\al_{i,j}=\v (Y_{i,j})=\v (Y_i-\sum_{k=1}^{j} c_{i,k}Y_1^{m_{i,k}})=\v (Y_{i,j-1}-c_{i,j}Y_1^{m_{i,j}}) >\v (Y_{i,j-1})=\al_{i,j-1}$$
Then we can obtain an infinite sequence of variables
$$
\begin{array}{rcl}
Y_{1,j} & = & Y_{1,j} \\
Y_{i,j} & = & Y_i-\sum_{k=1}^j c_{i,k}Y_1^{m_{i,k}},\
i=2,\ldots ,n,
\end{array}
$$
with $\v (Y_{i,j})>\v (Y_{i,j-1})$ for all $i,\ j$. So any sequence of partial sums of the series
$$Y_i-\sum_{k=1}^{\infty} c_{i,k}Y_1^{m_{i,k}},\ \forall
i=2,\ldots ,n$$
have strictly increasing values. Then these series converge to zero in the complete ring $\R$, so
$$Y_i=\sum_{k=1}^{\infty}
c_{i,k}Y_1^{m_{i,k}},\ \forall i=2,\ldots ,n.$$
Let $f(Y_1,\ldots ,Y_n)\in K_n$, then
$$v(f)=\v\( f\( Y_1,\sum_{k=1}^{\infty} c_{2,k}Y_1^{m_{2,k}}, \ldots ,\sum_{k=1}^{\infty} c_{n,k}Y_1^{m_{n,k}}\)\) =m\cdot
v(Y_1).$$
In this situation, the  group of $v$ is $v(Y_1)\cdot\ZZ$ (see \cite{Bri2}) but as the group is assumed to be $\ZZ$, $\v (Y_1)=1$. $\qed$



\section{Transcendental and algebraic elements of $\D$}

In the following sections we give a procedure to construct the residue field $\d$ of a discrete valuation of $K_n\vert k$, as a transcendental extension of $k$.

Before the describing procedure we have the following remark about the $k-$section $\si$.

\begin{rem}{\label{remark4}}
We are going to search all the variables for those residues which generate the extension $k\subset\D$. Hence we will have to move from $\D$ to $\R$ using the $k-$section $\si$ of the natural homomorphism $\R\to\D$.

As $\R$ is a complete ring, by Hensel's Lemma we know that $\D$ is a subfield of $\R$, in fact $\si (\D )$ is a subfield of $\R$ isomorphic to $\D$.

Since $\d$ and $\D$ are isomorphic, given $\om +\M\in\D$ we can find $f\in\r$ such that $f+\M =\om +\M$. So one may think that we can forget the $k$--section $\si$ because there exists a representative $f\in\r$ for each class $\om +\M$. This fact would mean that $\d$ is isomorphic of a subfield of $\r$ and this may not be true.

For example, the ring $R=\RR \lcor X\rcor_{(X^2+1)}$ is a rank-one discrete valuation ring, the completion is $\widehat{R}=\CC \lcor X\rcor$ and the residue field is $\CC\simeq R/(X^2+1)$, obviously $\CC\subset\widehat{R}$ but $\CC\not\subset R$.

So we need both the section and Hensel's Lemma to find $\d$ (or $\D$) as a subfield of $\R$ although we are able to find representatives in $\r$ for each element of $\d$.
\end{rem}

\begin{prop}{\label{prop4}}
Let us consider the diagram
\begin{center}
\unitlength=0.75mm
\linethickness{0.4pt}
\begin{picture}(49.00,59.67)
\vspace{.5cm} \put(10.00,4.00){\makebox(0,0)[cc]{$k$}}
\put(29.50,4.33){\vector(1,0){14.00}}
\put(29.50,4.33){\vector(-1,0){14.00}}
\put(28.00,6.33){\makebox(0,0)[cc]{{\footnotesize $id$}}}
\put(49.00,4.00){\makebox(0,0)[cc]{$k$}}
\put(10.00,7.00){\vector(0,1){15.00}}
\put(49.00,7.00){\vector(0,1){15.00}}
\put(10.00,27.67){\makebox(0,0)[cc]{$\FF$}}
\put(10.00,32.33){\vector(0,1){15.00}}
\put(49.00,32.33){\vector(0,1){15.00}}
\put(10.00,52.67){\makebox(0,0)[cc]{$\R$}}
\put(49.00,52.67){\makebox(0,0)[cc]{$\D$}}
\put(16.00,26.00){\vector(1,0){28.00}}
\put(44.00,28.00){\vector(-1,0){28.00}}
\put(28.00,24.00){\makebox(0,0)[cc]{{\footnotesize $\vp$}}}
\put(28.00,30.00){\makebox(0,0)[cc]{{\footnotesize $\si$}}}
\put(49.00,27.67){\makebox(0,0)[cc]{$\FF'$}}
\end{picture}
\end{center}
where $\FF$ and $\FF '$ are subfields of $\R$ and $\D$ respectively. If $\om +\M\in\D$, $\v (\om )=0$, is a transcendental (resp. algebraic) element over $\FF '$, there exists a $k-$section of $\vp$ which extends $\si$ and $\si (\om +\M )$ is transcendental (resp. algebraic) over $\FF$.
\end{prop}

\noindent {\it Proof.}

Let us suppose $\om +\M$ to be transcendental over $\FF '$. Let $\si\colon\D\to\R$ any section of $\vp$ extending $\si\colon\FF '\to \FF$. Let $f(X)\in\FF [X]$ be a non-zero polynomial. Let us put
$$f(X)=\sum_{i=0}^n\si (a_i')X^i,\ a_i'\in\FF '.$$
Then
$$f(\si
(\om +\M ))=\sum_{i=0}^n\si (a_i')\si (\om +\M
)^i=\si\(\sum_{i=0}^na_i'(\om +\M )^i\)\ne 0$$
because $\om +\M$ is transcendental over $\FF '$. So we have proved that $\si (\om + \M)$ is transcendental over $\FF$ if $\om +\M$ is transcendental over $\FF '$.

Now we suppose that $\al +\M\in\D$ is an algebraic element over $\FF '$, with $\v (\al )=0$ (i.e. $\al +\M\ne 0$). Let $$\fb (X) = X^m+\be_1X^{m-1}+\cdots +\be_m\in\FF '[X]$$
be its minimal polynomial over $\FF '$. Let us consider the polynomial
$$f(X) = X^m+b_1X^{m-1}+\cdots +b_m\in\FF[X],\mbox{ with } b_i=\si (\be_i).$$
By a corollary of Hensel's Lemma (\cite{ZSII}, corollary 1, page 279) we know that there exists $a\in\R$ such that $a$ is a simple root of $f(X)$ and $\varphi (a)=\al +\M$. As $\vp\si =id$, $f(X)$ is the minimal polynomial of $a$, so we can extend $\si :\FF '[\al +\M ]\to\FF [a]$. Then we have
\vspace{0.5cm}
\begin{center}
\unitlength=0.75mm
\linethickness{0.4pt}
\begin{picture}(49.00,77.34)
\vspace{.5cm} \put(10.00,4.00){\makebox(0,0)[cc]{$k$}}
\put(29.50,4.33){\vector(1,0){14.00}}
\put(29.50,4.33){\vector(-1,0){14.00}}
\put(28.00,6.33){\makebox(0,0)[cc]{{\footnotesize $id$}}}
\put(49.00,4.00){\makebox(0,0)[cc]{$k$}}
\put(10.00,7.00){\vector(0,1){15.00}}
\put(49.00,7.00){\vector(0,1){15.00}}
\put(10.00,27.67){\makebox(0,0)[cc]{$\FF$}}
\put(10.00,32.33){\vector(0,1){15.00}}
\put(49.00,32.33){\vector(0,1){15.00}}
\put(10.00,51.67){\makebox(0,0)[cc]{$\FF(a)$}}
\put(49.00,51.67){\makebox(0,0)[cc]{$\FF'(\al +\M )$}}
\put(10.00,57.00){\vector(0,1){15.00}}
\put(49.00,57.00){\vector(0,1){15.00}}
\put(10.00,77.34){\makebox(0,0)[cc]{$\R$}}
\put(49.00,77.34){\makebox(0,0)[cc]{$\D$}}
\put(16.00,26.00){\vector(1,0){28.00}}
\put(44.00,28.00){\vector(-1,0){28.00}}
\put(28.00,24.00){\makebox(0,0)[cc]{{\footnotesize $\vp$}}}
\put(28.00,30.00){\makebox(0,0)[cc]{{\footnotesize $\si$}}}
\put(49.00,27.67){\makebox(0,0)[cc]{$\FF'$}}
\put(20.00,50.67){\vector(1,0){15.00}}
\put(35.00,52.67){\vector(-1,0){15.00}}
\put(28.00,48.67){\makebox(0,0)[cc]{{\footnotesize $\vp$}}}
\put(28.00,54.67){\makebox(0,0)[cc]{{\footnotesize $\si$}}}
\end{picture}
\end{center}

Let us consider the set
$$\Om =\{ (\FF_1,\si_1)\vert\FF_1\supset\FF\mbox{ and }\si_1\mbox{ extends }\si\}$$
partially ordered by
$$(\FF_1,\si_1)<(\FF_2,\si_2)\iff\FF_1\subset\FF_2\mbox{ and  } \si_{2\vert\FF_1}=\si_1.$$
By Zorn's Lemma there exists a maximal element $(\LL ,\si ')\in\Om$, and again by another corollary of Hensel's Lemma (\cite{ZSII}, corollary 2, page 280) we have $\vp (\LL )=\D$. So we can extend $\si$ to a $k-$section $\si '$ of $\vp$ in such a way that $a=\si '(\al +\M)$ is an algebraic element over $\FF$.$\qed$


\section{A first transcendental residue.}

We devote this section to finding a first transcendental residue of $\d$ over $k$. Note that this preliminary transformations construct the residue field in the case $n=2$.

For a moment we shall suppose that the field $k$ is algebraically closed.

\begin{lem}{\label{lema6}}
There exists a finite number of monoidal transformations and changes of coordinates that constructs $n$ elements $Y_1,\ldots ,Y_n$ such that $v(Y_i)=v(Y_1)=\al$ and the residue $Y_2/Y_1+\M$ is not in $k$.
\end{lem}

\noindent{\it Proof.}
We can suppose that we have applied lemma \ref{lema1} to obtain $Y_1.\ldots ,Y_n$ such that $v(Y_i)=\al$ for all $i=1, \ldots ,n$.

In this situation $v(Y_i/Y_j)=0$, so $0\ne (Y_i/Y_j)+\m\in\d$. If this residue lies in $k$ then there exists $a_{i,1}\in k$ such that
$$\frac{Y_i}{Y_j}+\m = a_{i,1}+\m ,$$
so
$$\frac{Y_i}{Y_j}-a_{i,1}=\frac{Y_i-a_{i,1}Y_j}{Y_j}\in\m ,$$
and then
$$v\(\frac{Y_i-a_{i,1}Y_j}{Y_j}\) >0.$$
So we have $v(Y_i-a_{i,1}Y_j)=\al_1>\al$. If $\al$ divides to $\al_1$ then $\al_1=r_1\al$ with $r_1\geq 2$ and
$$v\(\frac{Y_i-a_{i,1}Y_j}{Y_j^{r_1}}\) =0.$$
If the residue of this element lies too in $k$, then exist $a_{i,r_1}\in k$ such that
$$v(Y_i-a_{i,1}Y_j-a_{i,r_1}Y_j^{r_1})=\al_2>\al_1.$$
If $\al$ divides to $\al_2$ then $\al_2=r_2\al$ with $r_2>r_1$ and we can repeat this operation.

The above procedure is finite for at least one pair $(i,j)$. We know (\cite{Bri2}) that any discrete valuation of $k((X_1,X_2))$ has dimension 1, so the restriction, $v'$, of our valuation $v$ to the field $k((X_1,X_2))$ is a valuation of dimension 1, and the dimension of $v$ is greater or equal than 1, because a transcendental residue of $v'$ over $k$ is a transcendental residue of $v$ too. If the procedure didn't finish for all $(i,j)$ then all the residues of $v$ would be in $k$, so the dimension of $v$ would be 0 and we have a contradiction. So we can suppose that the above procedure ends for $(1,2)$ by reordering the variables if necessary.

Hence there exists a first transcendental residue. We can apply the above procedure to the variables $Y_1,Y_2$, and so we have the transformations:
$$Z_i=Y_i,\ i\ne 2$$
$$Z_2= Y_2-\sum_{i=1}^{s_2}a_{2,i}Y_1^i,$$
such that one of the following two situations occurs:

a) $v(Y_1)$ divides $v(Z_2)$ and the residue of $Z_2/Y_1^r$ is not in $k$ with $v(Z_2)=r\cdot v(Y_1)$.

b) $v(Y_1)$ does not divide $v(Z_2)$.

In case a), we make the transformation
$$Z_2=Y_2-\sum_{i=1}^{s_2}a_{2,i}Y_1^i,$$
and apply lemma \ref{lema1} to obtain elements with the same values. We note these elements by $Y_1,\ldots ,Y_n$ again in order not to complicate the notation. So, after this procedure, we have a transcendental element $u_2=\si (Y_2/Y_1+\M)$ over $k$.

In case b) we make the same transformation and go back to the beginning of the proof.

Anyway this procedure stops, because the value of the variables are greater or equal than 1.

Then we can suppose that, after a finite number of transformations, we have $n$ elements $Y_1,\ldots ,Y_n$ such that $v(Y_i)=v(Y_1)=\al$ and the residue $Y_2/Y_1+\M$ is not in $k$.$\qed$

\begin{rem}{\label{notalgebraic}}
If $k$ is not algebraically closed the previous procedure gives us an algebraic extension $k\subset L=k[a_{j,i}]$, where $a_{j,i}$ are all the algebraic residues found in the procedure, and $n$ elements $Y_1,\ldots ,Y_n$ such that $\v (Y_i)=\v (Y_1)=\al$ and the residue $Y_2/Y_1+\M$ is transcendental over $L$.
\end{rem}

\begin{exmp}
Let $v=\nu_t\circ\Psi$ the discrete valuation of $\CC ((X_1,X_2))\vert\CC$ defined by the embedding
$$\begin{array}{rcl}
\Psi :\CC\lcor X_1,X_2\rcor & \lto & \CC (u)\lcor t\rcor \\
X_1 & \lmapsto & t\\
X_2 & \lmapsto & t+t^3+\sum_{i\geq 1}u^it^{i+3}
\end{array}
$$
with $u$ and $t$ independent variables over $\CC$.

The residue $X_2/X_1+\m =1+\m$, because $v(X_2-X_1)=3>1$. So we have
$$v\left(\frac{X_2-X_1}{X_1^3}\right) =0.$$
The residue
$$\frac{X_2-X_1}{X_1^3}+\m =1+\m$$
too, because $v(X_2-X_1-X_1^3)=4>3$. So we have
$$v\left(\frac{X_2-X_1-X_1^3}{X_1^4}\right) =0.$$
As $\Psi ((X_2-X_1-X_1^3)/X_1^4)=u$ and $u$ is trancendental over $\CC$, then
$$\frac{X_2-X_1-X_1^3}{X_1^4}+\m\notin \CC$$
and this is a first transcendental residue of $\d$ over $\CC$.

In this situation we can do the transformation
$$\begin{array}{rcl}
\CC\lcor X_1,X_2\rcor & \lto & \CC\lcor Y_1,Y_2\rcor\\
X_1 & \lmapsto & Y_1\\
X_2 & \lmapsto & Y_2Y_1^3+Y_1+Y_1^3
\end{array}
$$
to obtain elements $\{ Y_1,Y_2\}$ such that $\Psi (Y_1)=t$ and $\Psi (Y_2)=\sum_{i\geq 1}u^it^i$. So $v(Y_2)=v(Y_1)=1$ and the residue
$$\frac{Y_2}{Y_1}+\m =\frac{X_2-X_1-X_1^3}{X_1^4}+\m$$
is not in $\CC$.

In this example the extension of the valuation $v$ to the field $\CC ((Y_1,Y_2))$ is the usual order function. Theorem \ref{thm3} says that, for $n=2$, we always have this.
\end{exmp}

We end up the section with some specific arguments for the case $n=2$.

The proof of the following lemma is straightforward from (\cite{Bri2}, theorem 2.4):

\begin{lem}
Let $v$ be a discrete valuation of $K_n\vert k$. If $v$ is such that $v(f_r)=r\al$ for all forms $f_r$ of degree $r$ with respect to the usual degree, then the group of $v$ is $\al\cdot\ZZ$.
\end{lem}


\begin{thm}{\label{thm3}}
In the case $n=2$, the extension of the valuation $v$ to the field $k((Y_1,Y_2))$ is the usual order function.
\end{thm}

\noindent{\it Proof.}
After a finite number of transformations we are in the situation of the end of the previous proof. Obviously, if $n=2$,
$k((Y_1,Y_2))\subset\R$ so $v$ can be extended to a valuation $v'$ over $k((Y_1,Y_2))$ such that $\Delta_{v'}=\Delta_v=\Delta_{\v}$.
Let $\si :\D\to\R$ a $k-$section of $\R\to\D$, $u_2=\si (Y_2/Y_1+\M )$, $h\ne 0$ a form of degree $r$ and $\gamma =Y_2-u_2Y_1$. From the construction procedure of $u_2$ we know that $\v (\gamma )> \al$ (remember $\al =v'(Y_1)$). Then
$$h(Y_1,Y_2)=h(Y_1,u_2Y_1+\gamma
)=Y_1^rh(1,u_2)+\gamma ',$$
where $\gamma '$ is such that $v'(\gamma ')>r\al$. As $u_2\notin k$, $u_2$ is transcendental over $k$, so $h(1,u_2)\ne 0$ and $v'(h)=r\al$. By the previous lemma, the group of $v'$ is $\al\cdot\ZZ$, so $\al =1$ and $v'$ is the usual order function. $\qed$

\section{The general case}

Let us move to the general case. Assume $n>2$ and suppose we have applied the procedure of lemma \ref{lema6} in the transcendental case to find $Y_1,\ldots ,Y_n\in\K$ such that

a) The value of these elements are $\al\in\ZZ$.

b) The residue of $Y_2/Y_1$ is transcendental over $L$, where the extension $k\subset L$ is algebraic.

This section and the next one describe the transformations that we have to do in order to construct the residue field of $v$.

\begin{rem}
Let $\Delta_2=L(Y_2/Y_1+\M )$ a transcendental extension of $k$ of transcendence degree 1. Let $\si_2:\De_2\to L(Y_2/Y_1)$ defined by
$$\si_2\(\frac{Y_2}{Y_1}+\M\) =\frac{Y_2}{Y_1}=u_2.$$
We know that there exists a $k-$section $\si$ which extends $\si_2$ in the sense of proposition \ref{prop4}.
\end{rem}

\begin{rem}
Let  us suppose that the residue of $Y_3/Y_1$ is algebraic over $\Delta_2$, and let $u_{3,1}$ be its image by $\si$. Then $\v (Y_3-u_{3,1}Y_1)=\al_1>\al$. If $\al$ divides $\al_1$ then there exists $u_{3,r}\in\im (\si )$ and $r>1$ such that $\v (Y_3-u_{3,1}Y_1-u_{3,r}Y_1^r)=\al_2>\al_1$. Let us suppose that $u_{3,r}$ is algebraic over $\Delta_2$ too and $\al$ divides $\al_2$. Then we find ourselves in one of the three situations shown in the following items.

{\bf (Situation 1)} After a finite number of transformations, we obtain a value $\al_s$ such that it is not divided by $\al$. Then we make the transformation
$$Z_3=Y_3-\sum_{j=1}^su_{3,j}Y_1^j,$$
with $u_{3,j}$ algebraic over $\Delta_2$ for all $j=1,\ldots ,s$. So we have to apply transformations to find elements with the same values and begin with all the procedure described in this section. When this happens, the values of the elements decrease, so we can suppose that after a finite number of transformations we have reached a strictly minimal value. In fact this value should be 1, because we are assuming that the value group is $\ZZ$. We shall denote these elements by $Y_1,\ldots ,Y_n$ in order not to complicate the notation. So we can suppose that this situation will never happen again for any variable.

\noindent {\bf (Situation 2)} After a finite number of steps, we have a transcendental residue of $\Delta_2$. Let us denote this residue by $u_3$. This means
$$Z_3=Y_3-\sum_{j=1}^{s_3}u_{3,j}Y_1^j,$$
where the elements $\{ u_{3,j}\}_{j=1}^{s_3}$ are algebraic over $\Delta_2$ and $u_3=\si (Z_3/Y_1^{\v (Z_3)}+\M )$ is transcendental over $\Delta_2$. We shall note $\Delta_3 =L(u_2,\{ u_{3,j}\}_{j=1}^{s_3},u_3)$.

In this situation, if $n=3$ we can apply monoidal transformations to obtain elements with the same values. We will denote these elements again by $\{ Y_1,Y_2 ,Y_3\}$. The extension of the valuation $v$ to the field $\overline{L}((Y_1,Y_2,Y_3))$ with $\overline{L}=L(\{u_{3,j}\}_{j=1}^{s_3})$, is the usual order function, as in case $n=2$ (theorem \ref{thm3}).

\noindent {\bf (Situation 3)} All the residues obtained are algebraic elements. Then we take $\Delta_3=\Delta_2(\{ u_{3,j}\}_{j\geq 1})$, an algebraic extension of $\Delta_2$.
\end{rem}


\begin{rem}
Let us suppose that we have followed the previous construction with each element $Y_4,\ldots ,Y_{i-1}$, so we have a field
$$\De_{i-1}=L(u_2,\zeta_3,\ldots ,\zeta_{i-1})\subset \si (\D ),$$
where each $\zeta_k$ is:

- either $\{\{ u_{k,j}\}_{j=1}^{s_k},u_k\}$ if $\{ u_{k,j}\}_{j=1}^{s_k}$ are algebraic over $\De_{k-1}$ and $u_k=\si ((Z_k/Y_1^{\v (Z_k)})+\M )$ is a transcendental element over $\De_{k-1}$ (i.e. situation 2),

- or $\De_{k-1}\subset\De_{k-1}(\{ u_{k,j}\}_{j\geq 1})$ is an algebraic extension (i.e. situation 3).

So we have two possible situations concerning variable $Y_i$:

{\bf 1)} There exists a transformation
$$Z_i=Y_i-\sum_{j=1}^{s_i}u_{i,j}Y_1^j,$$
where the elements $u_{i,j}$ are algebraic over $\De_{i-1}$ and $u_i=\si ((Z_i/Y_1^{\v (Z_i)})+\M )$ is a transcendental element over $\De_{i-1}$. So we have the transcendental extension
$$\De_{i-1}\subset\De_{i-1}(\{u_{i,j}\}_{j=1}^{s_i},u_i)=\De_i.$$

{\bf 2)} All the elements $u_{i,j}$ we have constructed are algebraic over $\De_{i-1}$, so we have the algebraic extension
$$\De_{i-1}\subset\De_{i-1}(\{u_{i,j}\}_{j\geq 1})=\De_i.$$
\end{rem}

\begin{rem}
We have given a procedure to  construct elements $\{ Y_1,\ldots ,Y_n\}$ such that they satisfy these important properties:
\begin{enumerate}
\item After reordering if necessary, we can suppose that the first $m$ elements give us all the transcendental residues over $k$, i.e. the residue of each $Y_i/Y_1$ is transcendental over $\De_{i-1}$ with $i=2,\ldots ,m$. So the rest of variables $Y_{m+1},\ldots ,Y_n$ are such that we enter in situation {\bf 2)}.
\item With the usual notations, the extension
$$\De_m\subset\De_m\(\{u_{i,j}\}_{j\geq 1}\),\ i=m+1,\ldots ,n$$
is algebraic.
\item We are assuming that the dimension of $v$, $\dim (v)$, is equal to $m-1$.
\end{enumerate}
\end{rem}


\begin{thm}{\label{theorem9}}
The residue field of $v$ as subfield of $\R$ is
$$\De_n=L\( u_2,\{u_{3,j}\}_{j=1}^{s_3},u_3,\ldots
,\{u_{m,j}\}_{j=1}^{s_m},u_m\) \(\{u_{m+1,j}\}_{j\geq 1},\ldots
,\{u_{n,j}\}_{j\geq 1}\) ,$$
and the dimension of $v$, i.e. the transcendence degree of $\De_n$ over $k$, is $m-1$.
\end{thm}

\noindent{\it Proof.}
In this section we have given a construction by writing the elements $Y_i$ depending on $Y_1$ and some transcendental and algebraic residues. So we have constructed a map
$$
\begin{array}{rcl}
\vp ':L\lcor \Y \rcor & \longrightarrow & \Delta_n\lcor t \rcor \\
Y_1 & \lmapsto & t \\
Y_i & \lmapsto & u_it,\ i=2,\ldots ,m\\
Y_k & \lmapsto & \sum_{j\geq 1}u_{k,j}t^j,\ u_{k,1}\ne 0,\
k=m+1,\ldots ,n.
\end{array}
$$
This map is not injective in the general case, but we know that $v=\nu_t\circ\phi '_{\vert K_n}$, where $\phi '$ is the extension of $\vp '$ to the quotient field. So the residue field of $v$ is equal to the residue field of $\nu_t$, i.e. $\De_n$. $\qed$

A straightforward consequence of this theorem is the following well-known result

\begin{cor}{\label{cor8}}
The usual order function over $K_n$ has dimension $n-1$, i.e. the transcendence degree of its residue field over $k$ is $n-1$.
\end{cor}

\noindent{\it Proof.}
Let $\nu$ be the usual order function over $K_n$. All the residues $X_i/X_1+{\got m}_{\nu}$ are transcendental over $k(X_2/X_1+{\got m}_{\nu},\ldots , X_{i-1}/X_1+{\got m}_{\nu})$: if this were not the case, there would exist $u_i\in\si (\Delta_{\nu})$ such that $\nu (X_i-u_iX_1)>1$ and $\nu$ would not be an order function. So $\Delta_{\nu}=k(X_2/X_1,\ldots ,X_{n}/X_1)$. $\qed$

\section{Explicit construction of the residue field: an example}

In order to compute explicitly the residue field of a valuation we need to construct a section $\si :\D\to\R$ as in proposition \ref{prop4}. This procedure is not constructive in general. As in section 1, if the valuation is given as a composition $v=\nu_t\circ\Psi$, where $\Psi :k\lcor\X\rcor\to\Delta\lcor t\rcor$ is an injective homomorphism and $\nu_t$ is the order funcion in $\Delta\lcor t\rcor$, then we can construct $\si$ using the coefficients $a_{i,j}\in\Delta$ of $\Psi (X_i)=\sum_{j\geq 1}a_{i,j}t^j$.

(Of course, explicit does not mean effective because we are working with series $\sum_{j\geq 1}a_{i,j}t^j$ and this input is not finite).

\begin{exmp}{\label{ejemplo}}
Let us consider the embedding
$$
\begin{array}{rcl}
\Psi :  \CC\lcor X_1,X_2,X_3,X_4,X_5\rcor & \lto &
        \Delta\lcor t\rcor \\
X_1 & \lmapsto & t \\
X_2 & \lmapsto & T_2t \\
X_3 & \lmapsto & T_2^2t+T_2t^2+T_3t^3 \\
X_4 & \lmapsto & T_2^3t+T_2^2t^2+T_3t^3+T_4t^4\\
X_5 & \lmapsto & T_2t\sum_{j\geq 1} (T_4^{1/p}t)^j,
\end{array}
$$
with $t$, $T_2$, $T_3$ and $T_4$ variables over $\CC$, $p\in\ZZ$ prime and $\Delta$ is a field such that $\overline{\CC (T_4)}(T_2,T_3)\subseteq\Delta$. $\overline{\CC (T_4)}$ is the algebraic closure of $\CC (T_4)$. We are going to denote  its extension to the quotient fields by $\Psi$. The composition of this injective homomorphism with the order function in $t$ gives a discrete valuation of $\CC ((X_1,X_2,X_3,X_4,X_5))\vert\CC$, $v=\nu_t\circ\Psi$. The residues of $X_i/X_1$ are not in $\CC$ for $i=2,3,4,5$.

Let us put $u_2=\si (X_2/X_1+\M )$, a transcendental element over $\CC$. By proposition \ref{prop4} we know how to construct $\si$ step by step, so let take us $u_2=X_2/X_1$ and $\De_2=\CC(u_2)$.

The residue $X_3/X_1+\M$ is algebraic over $\CC (u_2)$, in fact
$$\frac{X_3}{X_1}+\M =\frac{X_2^2}{X_1^2}+\M .$$
So we can take $u_{3,1}=\si ((X_3/X_1)+\M )=u_2^2$. The value of $X_3-u_{3,1}X_1$ is 2, therefore we have to see if the residue
$$\frac{X_3-u_{3,1}X_1}{X_1^2}+\M$$
is algebraic over $\CC (u_2)$. We have that $$\frac{X_3-u_{3,1}X_1}{X_1^2}+\M = \frac{X_2}{X_1}+\M
,$$
so it is algebraic and we can take $u_{3,2}=u_2$. Now $v(X_3-u_{3,1}X_1-u_{3,2}X_1^2)=3$ and we have to check if
$$\frac{X_3-u_{3,1}X_1-u_{3,2}X_1^2}{X_1^3}+\M$$
is algebraic over $\De_2$. In this case, as
$$\Psi\(\frac{X_3-u_{3,1}X_1-u_{3,2}X_1^2}{X_1^3}+\M\)= T_3 ,$$
this residue is transcendental. So we take
$$u_3=\si\(\frac{X_3-u_{3,1}X_1-u_{3,2}X_1^2}{X_1^3}+\M\)
=\frac{X_1X_3-X_2^2-X_1^2X_2}{X_1^4} .$$
Let us take $\De_3=\CC (u_2,u_3)$.

We have to apply this procedure to $X_4$. The residue $X_4/X_1+\M$ is algebraic over $\De_3$ because
$$\frac{X_4}{X_1}+\M =\frac{X_2^3}{X_1^3}+\M ,$$
so we can take $u_{4,1}=\si ((X_4/X_1)+\M )=u_2^3\in\De_3.$ Now $v(X_4-u_{4,1}X_1)=2$, and we have to check what happens with the residue
$$\frac{X_4-u_{4,1}X_1}{X_1^2}+\M .$$
As
$$\frac{X_4-u_{4,1}X_1}{X_1^2}+\M = \frac{X_1^2}{X_2^2}+\M ,$$
it holds
$$u_{4,2}=\si\(\frac{X_4-u_{4,1}X_1}{X_1^2}+\M\) =u_2^2.$$
Clearly $v(X_4-u_{4,1}X_1-u_{4,2}X_1^2)=3$ and
$$\frac{X_4-u_{4,1}X_1-u_{4,2}X_1^2}{X_1^3}+\M =\frac{X_1X_3-X_2^2-X_1^2X_2}{X_1^4}+\M ,$$
therefore
$$u_{4,3}=\si\(\frac{X_4-u_{4,1}X_1-u_{4,2}X_1^2}{X_1^3}+\M\) =u_3.$$
The following residue is transcendental because
$$v(X_4-u_{4,1}X_1-u_{4,2}X_1^2-u_{4,3}X_1^3)=4$$
and
$$\Psi\(\frac{X_4-u_{4,1}X_1-u_{4,2}X_1^2-u_{4,3}X_1^3}{X_1^4}\)=T_4.$$
Then we can take
$$u_4=\si\(\frac{X_4-u_{4,1}X_1-u_{4,2}X_1^2-u_{4,3}X_1^3}{X_1^4}+\M\) =$$
$$=\frac{X_1^2X_4-X_3^2-X_1^2X_2^2-X_1^2X_3-X_1X_2^2-X_1^2X_2}{X_1^6}.$$
So $\De_4=\CC (u_2,u_3,u_4)$.

With the variable $X_5$ we obtain the next algebraic residues
$$u_{5,j}=\si\(\frac{X_5-u_{5,1}X_1-\cdots -u_{5,j-1}X_1^{j-1}}{X_1^j}+\M\) =u_4^{\frac{1}{p^j}}$$
for all $j\geq 1$. So we have $\Delta_5=\CC (u_2,u_3,u_4)(\{ u_4^{1/p^j}\}_{j\geq 1} )$, an algebraic extension of $\Delta_4$.

Then the residue field of $v$ is the subfield of $\R$
$$\d =\CC\(\frac{X_2}{X_1}+\m ,\frac{X_1X_3-X_2^2-X_1^2X_2}{X_1^4}+\m , \right.$$
$$\left.\frac{X_1^2X_4-X_3^2-X_1^2X_2^2-X_1^2X_3-X_1X_2^2-X_1^2X_2}{X_1^6}+\m\) \(\left\{\(\frac{X_2}{X_1}+\m\)^{\frac{1}{p^j}}\right\}_{j\geq 1}\) .$$
So we can find elements of $\r$ as representatives of the generators of the residual field, but $\d\not\subset\r$.

In this case, by the transformation
$$
\begin{array}{rcl}
X_1 & \lto & Y_1 \\
X_2 & \lto & Y_2 \\
X_3 & \lto & Y_1^2Y_3+u_{3,1}Y_1+u_{3,2}Y_1^2 \\
X_4 & \lto & Y_1^3Y_4+u_{4,1}Y_1+u_{4,2}Y_1^2+u_{4,3}Y_1^3 \\
X_5 & \lto & Y_5,
\end{array}
$$
we can extend the valuation $v$ to a discrete valuation $v'=\nu_t\Psi'$ of the field $\CC ((Y_1,Y_2,Y_3,Y_4,Y_5))$, with the injective homomorphism
$$
\begin{array}{rcl}
\Psi' :  \CC\lcor Y_1,Y_2,Y_3,Y_4,Y_5\rcor & \lto &
        \Delta\lcor t\rcor \\
Y_1 & \lmapsto & t \\
Y_i & \lmapsto & T_it,\ i=2,3,4\\
Y_5 & \lmapsto & \sum_{j\geq 1} (T_4^{1/p}t)^j.
\end{array}
$$
The restriction $v'_{\vert\CC ((Y_1,Y_2,Y_3,Y_4))}$ is the usual order function. This is not the general case because $\Psi'$ may not be injective.
\end{exmp}

\section{Rank one discrete valuations and order functions}

We can summarize the constructions of previous sections in the following theorem wich generalize the results of \cite{Bri2,Br-He}

\begin{thm}{\label{teorema212}}
Let $v$ be a discrete valuation of $K_n\vert k$, then
\begin{enumerate}
\item If the dimension of $v$ is $n-1$, we can embed $k\lcor\X\rcor$ into a ring $L\lcor\Y\rcor$, where $L\subset\si (\D )$ and the extended valuation of $v$ over the field $L((\Y ))$ is the usual order function.
\item If the dimension of $v$ is $m-1<n-1$, we can embed $k\lcor\X\rcor$ into a ring $L\lcor\Y\rcor$, where $L\subset\si (\D )$ and the restriction into $L((\Ym ))$ of the ``extended valuation'' of $v$ over $L((\Y ))$ is the usual order function.
\end{enumerate}
\end{thm}

\noindent{\it Proof.} We have the following map:
$$
\begin{array}{rcl}
\vp ':L\lcor \Y \rcor & \longrightarrow & \Delta_n\lcor t\rcor \\
Y_1 & \lmapsto & t \\
Y_i & \lmapsto & u_it,\ i=2,\ldots ,m \\
Y_k & \lmapsto & \sum_{j\geq 1}u_{k,j}t^j,\ u_{k,1}\ne 0,\
k=m+1,\ldots ,n,
\end{array}
$$
where $m-1$ is the dimension of $v$. Let $\phi '$ the extension of $\vp '$ to the quotient field. Let us prove the theorem:

\begin{enumerate}

\item In the case $m=n$, $\vp '(Y_i)=u_it$ for all $i=2,\ldots ,n$. Let $\nu_t$ be the usual order funtion over $\Delta_n ((t))$. The homomorphism $\vp '$ is injective and the valuation $v'=\nu_t\circ\phi '$ of $L((\Y ))$ is the usual order fuction over this field. Obviously $v'$ extends $v$.

\item If $m<n$ we can consider the elements $W_k=Y_k-\sum_{j\geq 1}u_{k,j}Y_1^j$. Hence we have $L ((\Y)) = L((Y_1,\ldots ,Y_m,W_{m+1},\ldots ,W_n))$. We define the discrete valuation of rank $n-\dim (v)=n-m+1$ over $L(( \Y ))$:
$$v' (Y_1)=\ldots
=v'(Y_m)= (0,\ldots ,0,1),$$
$$v'(W_{m+1})=(0,\ldots ,1,0),\ldots
,v'(W_n)=(1,0,\ldots ,0).$$
The restriction of this valuation to $K_n$ is a rank one discrete valuation, because the value of any element is in $0\times\cdots\times 0\times\ZZ$. In fact $v'(f)=(0,\ldots ,0,v(f))$ for all $f\in K_n$, so $v'$ ``extends'' $v$ in this sense. Obviously $v'_{\vert L((\Ym ))}$ is the usual order function. $\qed$
\end{enumerate}

\begin{rem}
Note that $(W_{m+1},\ldots ,W_{n})$ is the {\em implicit ideal of $v$} that appears in some works of M. Spivakovsky \cite{Sp}, B. Teissier \cite{Te} and both authors with F.J. Herrera and M.A. Olalla \cite{HOST}. This implicit ideal appears when we complete the ring after a monoidal transformation.
\end{rem}

For the case of valuations of dimension $n-1$, we can combine corollary \ref{cor8} and assertion 1 of the previous theorem:

\begin{cor}{\label{theorem8}}
Let $v$ be a discrete valuation of $K_n\vert k$. The following conditions are equivalent:

1) The transcendence degree of $\D$ over $k$ is $n-1$ (i.e. $\dim (v)=n-1$).

2) There exists a finite sequence of monoidal transformations and coordinates changes which take $v$ into an order function.
\end{cor}

\begin{exmp}
Let us consider the homomorphism
$$
\begin{array}{rcl}
\Psi :  \CC\lcor X_1,X_2,X_3,X_4,X_5\rcor & \lto &
        \Delta\lcor t\rcor \\
X_1 & \lmapsto & t \\
X_2 & \lmapsto & T_2t \\
X_3 & \lmapsto & T_2^2t+T_2t^2+T_3t^3 \\
X_4 & \lmapsto & T_2^3t+T_2^2t^2+T_3t^3+T_4t^4\\
X_5 & \lmapsto & T_2t\left(\sum_{j\geq 1} a_j(T_4t)^j\right) ,
\end{array}
$$
with $a_j\in\CC$ such that $\Psi$ is injective (we can take $\sum_{j\geq 1} a_j(T_4t)^j= e^{T_4t}-1$). Then the residue field of this valuation (see example \ref{ejemplo}) is
$$\d =\CC\(\frac{X_2}{X_1}+\m ,\frac{X_1X_3-X_2^2-X_1^2X_2}{X_1^4}+\m , \right.$$
$$\left.\frac{X_1^2X_4-X_3^2-X_1^2X_2^2-X_1^2X_3-X_1X_2^2-X_1^2X_2}{X_1^6}+\m\)\subset\R .$$
By the transformation (see example \ref{ejemplo})
$$
\begin{array}{rcl}
X_1 & \lto & Y_1 \\
X_2 & \lto & Y_2 \\
X_3 & \lto & Y_1^2Y_3+u_{3,1}Y_1+u_{3,2}Y_1^2 \\
X_4 & \lto & Y_1^3Y_4+u_{4,1}Y_1+u_{4,2}Y_1^2+u_{4,3}Y_1^3 \\
X_5 & \lto & Y_2Y_5,
\end{array}
$$
we obtain a new field $\CC ((Y_1,Y_2,Y_3,Y_4,Y_5))$, but we can not extend $v$ to this field because the homomorphism
$$
\begin{array}{rcl}
\Psi' :  \CC\lcor Y_1,Y_2,Y_3,Y_4,Y_5\rcor & \lto &
        \Delta\lcor t\rcor \\
Y_1 & \lmapsto & t \\
Y_i & \lmapsto & T_it,\ i=2,3,4 \\
Y_5 & \lmapsto & \sum_{j\geq 1} a_j(T_4t)^j
\end{array}
$$
is not injective. Then let us take $W_5=Y_5-\sum_{j\geq 1} a_j(Y_4)^j$ (because we can consider $T_4Y_1=Y_4$). Then $\CC ((Y_1,Y_2,Y_3,Y_4,Y_5))=\CC ((Y_1,Y_2,Y_3,Y_4,W_5))$ and the discrete valuation of rank 2 defined by $v'(Y_i)=(0,1)$ for $i=1,\ldots ,4$ and $v'(W_5)=(1,0)$ is such that for all $f\in\CC ((X_1,X_2,X_3,X_4,X_5))$ we have $v'(f)=(0,v(f))$ and $v'_{\vert\CC ((Y_1,Y_2,Y_3,Y_4))}$ is the usual order function.
\end{exmp}


\providecommand{\bysame}{\leavevmode\hbox to3em{\hrulefill}\thinspace}

\end{document}